\documentclass[11pt]{amsart}
\usepackage{amsmath}\usepackage{amssymb}
\usepackage{amsfonts}\usepackage{amscd}
\newcommand\op{\operatorname}

\newcommand\coker{\op{Coker}}

\newcommand\End{\op{End}}

\newcommand\GL{\op{GL}}

\newcommand\Hom{\op{Hom}}

\newcommand\image{\op{Im}}

\newcommand\Ker{\op{Ker}}

\newcommand\Proj{\op{Proj}}

\newcommand\quot{/\hskip -0.1cm /}
\newcommand\rank{\op{rank}}

\newcommand\SL{\op{SL}}

\newcommand\Spec{\op{Spec}\,}

\newcommand\varep{\varepsilon}

\newcommand{\bC}{{\bold C}}

\newcommand{\bG}{{\bold G}}

\newcommand{\bL}{{\bold L}}

\newcommand{\bP}{{\bold P}}
\newcommand{\bQ}{{\bold Q}}

\newcommand{\bZ}{{\bold Z}}

\newtheorem{thm}[subsection]{Theorem}

\newtheorem{lemma}[subsection]{Lemma}
\newtheorem{cor}[subsection]{Corollary}

\newtheorem{conjecture}[subsection]{Conjecture}

\theoremstyle{definition}
\newtheorem{defn}[subsection]{Definition}

\newtheorem{rem}[subsection]{Remark}

\ifx\undefined\pf
\newcommand{\pf}{\proof}
\fi
\ifx\undefined\cal
\newcommand{\cal}{\mathcal}
\fi
\ifx\undefined\Bbb
\newcommand{\Bbb}{\mathbb}
\fi
\ifx\undefined\bold
\newcommand{\bold}{\mathbf}
\fi
\ifx\undefined\frak
\newcommand{\frak}{\mathfrak}
\fi
\begin{document}
\setcounter{tocdepth}{1}
\setcounter{section}{0}
\title[Conics on a generic hypersurface]
{Conics on a generic hypersurface}
\author{Masao Jinzenji, Iku Nakamura and Yasuki Suzuki}
\maketitle
\begin{abstract}
In this paper, we compute the contributions from double cover maps to genus 0 
degree 2 Gromov-Witten invariants of general type projective hypersurfaces. 
Our results correspond to a generalization of
Aspinwall-Morrison formula to general type hypersurfaces 
in some special cases.
\\
MSC-class: 14H99, 14N35, 32G20 
\end{abstract}
\section{Introduction}
\label{sec:Introduction}In this paper, we discuss a generalization
of the multiple cover formula for rational
Gromov-Witten invariants of Calabi-Yau manifolds \cite{AM}, \cite{manin} 
to double cover maps of a line $L$ on a degree $k$
hypersurface $M_{N}^{k}$ in $\bP^{N-1}$.  N\"aively,
for a given finite set of elements
$\alpha_{j}\in H^{*}(M_{N}^{k},\bZ)$,
the rational Gromov-Witten invariant $\langle {\cal O}_{\alpha_{1}}{\cal
O}_{\alpha_{2}}\cdots
{\cal O}_{\alpha_{n}}\rangle_{0,d}$ of $M_{N}^{k}$
counts the number of degree $d$ (possibly singular and reducible) rational
curves on $M_{N}^{k}$  that
intersect real sub-manifolds of $M_{N}^{k}$ that are Poincar\'e-dual to
$\alpha_{j}$.

Recently, the mirror computation
of rational Gromov-Witten invariants of $M_{N}^{k}$
with negative first chern
class ($k-N>0$) was established in \cite{gc}, \cite{iri}, \cite{jin}.
Using the method presented in these articles, we can compute 
$\langle {\cal O}_{e^{m_{1}}}{\cal
O}_{e^{m_{2}}}\cdots
{\cal O}_{e^{m_{n}}}\rangle_{0,d}$
where $e$ is the generator of $H^{1,1}(M_{N}^{k},\bZ)$.
Briefly, mirror computation of $M_{N}^{k}\;\;(k>N)$ in \cite{jin}
goes as follows.
We start from the following ODE:
\begin{equation}
\biggl((\partial_{x})^{N-1}-k\cdot \exp(x)\cdot (k\partial_{x}+k-1)
(k\partial_{x}+k-2)
\cdots(k\partial_{x}+1)\biggr)w(x)=0,
\label{ifun}
\end{equation}
and construct the virtual Gauss-Manin system associated with (\ref{ifun}):
\begin{eqnarray}
\partial_{x}\tilde{\psi}_{N-2-m}(x)&=&\tilde{\psi}_{N-1-m}(x)+
\sum_{d=1}^{\infty}
\exp(dx)\cdot \tilde{L}_{m}^{N,k,d}\cdot\tilde{\psi}_{N-1-m-(N-k)d}(x),
\label{gm1}
\end{eqnarray}
where $m$ runs through all the integers and 
$\tilde{L}_{m}^{N,k,d}$ is non-zero 
only if $0\leq m\leq N-1+(k-N)d$. From the compatibility of (\ref{ifun})
and (\ref{gm1}), we can derive the recursive formulas that determine 
all the $\tilde{L}_{m}^{N,k,d}$'s: 
\begin{eqnarray*}
&&\sum_{n=0}^{k-1}\tilde{L}_{n}^{N,k,1}w^{n}=
k\cdot\prod_{j=1}^{k-1}(jw+(k-j)), \nonumber\\
&&\sum_{m=0}^{N-1+(k-N)d}\tilde{L}_{m}^{N,k,d}z^{m}=
\sum_{l=2}^{d}(-1)^{l}\sum_{0=i_{0}<\cdots<i_{l}=d}\times
\nonumber\\
&&\times\sum_{j_{l}=0}^{N-1+(k-N)d}
\cdots\sum_{j_{2}=0}^{j_{3}}\sum_{j_{1}=0}^{j_{2}}
\prod_{n=1}^{l}\biggl((\frac{i_{n-1}+(d-i_{n-1})z}{d})^{j_{n}-j_{n-1}}
\cdot \tilde{L}_{j_{n}+(N-k)i_{n-1}}^{N,k,i_{n}-i_{n-1}}\biggr).
\nonumber\\
\end{eqnarray*}
With these data, we can construct the formulas that represent rational 
three point 
Gromov-Witten invariant $\langle{\cal O}_{e}{\cal O}_{e^{N-2-m}}
{\cal O}_{e^{m-1-(k-N)d}}\rangle_{d}$ in terms of $\tilde{L}_{m}^{N,k,d}$.
These three point Gromov-Witten invariants are enough for reconstruction 
of all the rational Gromov-Witten invariants  
$\langle {\cal O}_{e^{m_{1}}}{\cal
O}_{e^{m_{2}}}\cdots
{\cal O}_{e^{m_{n}}}\rangle_{0,d}$ \cite{km}.
In particular, we obtain the following formula in the $d=2$ case:   
\begin{eqnarray}
&&\langle{\cal O}_{e}{\cal O}_{e^{N-2-m}}
{\cal O}_{e^{m-1-(k-N)2}}\rangle_{2}=\\
&&k\cdot\biggl(\tilde{L}_{n}^{N,k,2}-\tilde{L}_{1+2(k-N)}^{N,k,2}
-2\tilde{L}_{1+(k-N)}^{N,k,1}(\sum_{j=0}^{k-N}(\tilde{L}_{n-j}^{N,k,1} 
- \tilde{L}_{1+2(k-N)-j}^{N,k,1}))\biggr).\nonumber
\label{con}
\end{eqnarray} 
According
to the results of this procedure, rational three point 
Gromov-Witten invariants can be rational
numbers with large denominator if $k>N$ ,
in contrast to the Calabi-Yau case where rational
three point Gromov-Witten invariants are always integers.  \par
One of the reasons of this rationality (non-integrality) comes 
from the contributions of multiple cover maps to Gromov-Witten invariants.
In the Calabi-Yau case ($N=k$), for any divisor $m$ of $d$
there are some contributions from degree $m$ multiple cover
maps $\phi$ of a rational curve $\bP^{1}$ onto
a degree $\frac{d}{m}$ rational curve
$C\hookrightarrow M_{k}^{k}$.
The contributions from the multiple cover  maps
are expressed in terms of the virtual fundamental class of Gromov-Witten
invariants. 
Let $C$ be a  general degree $d$ rational curve in $M_{k}^{k}$.
Its normal bundle $N_{C/M_{k}^{k}}$ is decomposed into a direct sum of line
bundles as follows:
\begin{equation*}
N_{C/M_{k}^{k}}\simeq O_{C}(-1)\oplus O_{C}(-1)\oplus O_{C}^{\oplus
{(k-5)}}.
\end{equation*}

Let $\phi:\bP^{1}\to C$ be a holomorphic map of degree $m$.
Since the pull-back $\phi^{*}(N_{C/M_{k}^{k}})$ is given by
\begin{equation*}
\phi^{*}(N_{C/M_{k}^{k}})
\simeq O_{\bP^{1}}(-m)\oplus O_{\bP^{1}}(-m)\oplus O_{\bP^{1}}^{\oplus
{k-5}},
\end{equation*}    
 we obtain $h^{1}(\phi^{*}(N_{C/M_{k}^{k}}))=2m-2$. On the other hand, let
$\overline{M}_{0,0}(M,d)$ be the moduli space of $0$-pointed stable
maps of degree $d$ from genus $0$ curve to $M$. Then  the moduli space of
$\phi$ is 
the fiber space $\pi: \overline{M}_{0,0}(C,m)\rightarrow
\overline{M}_{0,0}(M_{k}^{k},\frac{d}{m})$, whose fibre
$\overline{M}_{0,0}(C,m)$ over $C$ (fixed) has
complex dimension $2m-2$. Then
the push-forward of the virtual fundamental class $\pi_{*}(c_{top}(H^{1}
(\phi^{*}N_{C/M_{k}^{k}})))$ can be computed only by intersection theory
on 
the fiber $\overline{M}_{0,0}(C,m)$, which turns out to be equal to
$\frac{1}{d^{3}}$.
This depends on neither the structure of the base
$\overline{M}_{0,0}(M_{k}^{k},
\frac{d}{m})$ nor the global structure of the fibration.

But when $k<N$, the situation is more complicated than $M_{k}^{k}$
because of negative first Chern class. Let us concentrate on the case of
$d=2, m=2$ that we discuss in this paper. In this case, $C$ is just a line
$L$ on the hypersurface $M_{N}^{k}$.  The moduli space
$\overline{M}_{0,0}(M_{N}^{k},1)$
is a sub-manifold of $\overline{M}_{0,0}(\bP^{N-1},1)$, while
$\overline{M}_{0,0}(\bP^{N-1},1)$
is the Grassmannian
$G(2, N)$, the moduli space of rank $2$
quotients  of $V=\bC^{N}$.
As will be shown later, for a generic line $L$, $N_{L/M_{N}^{k}}$ is
decomposed into 
\begin{equation*}
N_{L/M_{N}^{k}}\simeq O_{L}(-1)^{\oplus k-N+2}\oplus O_{L}^{\oplus
{2N-k-5}}.
\end{equation*}
By pulling back it by the degree $2$ map $\phi:\bP^{1}\rightarrow L$,
we obtain,
\begin{equation*}
\phi^{*}N_{L/M_{N}^{k}}\simeq O_{\bP^{1}}(-2)^{\oplus k-N+2}\oplus
O_{\bP^{1}}^{\oplus {2N-k-5}}.
\end{equation*} 
Therefore, $h^{1}(\phi^{*}(N_{L/M_{N}^{k}}))=k-N+2$, which is strictly
grater than two, 
the complex dimension of the fiber
$\overline{M}_{0,0}(L,2)$. Thus we need to know the global structure
of the fibration $\pi$
in order to compute the multiple cover contribution to degree $2$
rational Gromov-Witten invariants of $M_{N}^{k}$.

In order to estimate the contributions from double cover maps 
$\phi: \bP^{1}\rightarrow L$ to  
$\langle{\cal O}_{e^{a}}{\cal O}_{e^{b}}{\cal O}_{e^{c}}\rangle_{0,2}$, 
we first
computed the number of conics, that intersect cycles Poancar\'e dual to 
$e^{a}, e^{b}$ and $e^{c}$, on
$M_{N}^{k}$ 
(whose normal bundle are of the same type) by using the method in 
\cite{Katz2}. Then we  found the following formula by comparing these 
integers with the results obtained from (3):
\begin{eqnarray}
&&\langle{\cal O}_{e^{a}}{\cal O}_{e^{b}}{\cal O}_{e^{c}}\rangle_{0,2}
=(\mbox{number of corresponding conics})+\\
&&\int_{G(2,N)}c_{top}(S^{k}Q)\wedge\bigl[\frac{c(S^{k-1}Q)}
{1-\frac{1}{2}c_{1}(Q)}\bigr]_{k-N}\wedge\sigma_{a-1}\wedge\sigma_{b-1}
\wedge
\sigma_{c-1},\nonumber
\label{dc}
\end{eqnarray}
where $Q$ is the universal rank $2$ quotient bundle of $G(2,N)$,
$\sigma_{a}$ 
is a Schubert cycle defined by $\sum_{a={0}}^{\infty}\sigma_{a}:=
\frac{1}{c(Q^{\vee})}$ and $[*]_{k-N}$ is the operation of picking up degree
$2(k-N)$ part of Chern classes. \par
On the other hand, we have the following formula which directly follows 
from the definition of the virtual fundamental class of 
$\overline{M}_{0,0}(M_{N}^{k},2)$: 
\begin{eqnarray}
&&\langle{\cal O}_{e^{a}}{\cal O}_{e^{b}}{\cal O}_{e^{c}}\rangle_{0,2}
=(\mbox{number of corresponding conics})+\\
&&
8\int_{G(2,N)}c_{top}(S^{k}Q)\wedge\bigl[\pi_{*}
(c_{top}(H^{1}(\phi^{*}N_{L/M_{N}^{k}})))
\bigr]_{k-N}\wedge\sigma_{a-1}\wedge\sigma_{b-1}\wedge
\sigma_{c-1}.\nonumber
\label{dec}
\end{eqnarray}
where $\pi:\overline{M}_{0,0}(L,2)\rightarrow
\overline{M}_{0,0}(M_{N}^{k},1)$ is the natural projection.
Here, the factor $8$ comes from the divisor axiom of Gromov-Witten 
invariants.

In this paper, we prove the following formula
\begin{equation}
\pi_{*}(c_{top}(H^{1}
(\phi^{*}N_{L/M_{N}^{k}})))=\frac{1}{8}\bigl[\frac{c(S^{k-1}Q)}
{1-\frac{1}{2}c_{1}(Q)}\bigr]_{k-N}.
\label{re}
\end{equation}
By combining (5) with (\ref{re}), we can derive the formula 
(4) immediately. \par
 From (4), we see that
$\langle{\cal O}_{e^{a}}{\cal O}_{e^{b}}{\cal O}_{e^{c}}\rangle_{0,2}$ of
$M_{N}^{k}$ is a rational number with denominator at most $2^{k-N}$.
Therefore 
rationality (non-integrality) of the Gromov-Witten invariant
$\langle{\cal O}_{e^{a}}{\cal O}_{e^{b}}{\cal O}_{e^{c}}\rangle_{0,2}$ is
caused by the 
effect of multiple cover map in this case.
\par
We note here that the total
moduli space of double cover maps of lines is isomorphic to ${\bP}(S^{2}Q)$
over $G:=\overline{M}_{0,0}(M_{N}^{k},1)\hookrightarrow G(2,N)$, 
which is an algebraic $\bQ$-stack ${\bP}(S^{2}Q)^{stack}$
(in the sense of Mumford). As a consequence, the union of all
$H^{1}(\phi^{*}N_{C/M_{N}^{k}})$
turns out to be a coherent sheaf 
on ${\bP}(S^{2}Q)^{stack}$ with fractional
Chern class in (\ref{re}), 
as was suggested in \cite{thber}.  See \cite[Section 9]{Vistoli89}.

We also did some numerical experiments on degree $3$ Gromov-Witten
invariants 
of $M_{N}^{k}$ by using the results of \cite{ES}.
For  $k-N>0$, there is a new contribution from multiple cover maps to nodal
conics in $M_{N}^{k}$ that did not appear in the Calabi-Yau case. Therefore,
multiple cover map contributions are far more complicated than Calabi-Yau,
and we leave general analysis on this problem to future works.

This paper is organized as follows. In Section 1, we analyze characteristics
of moduli space of lines in $M_{N}^{k}$ and
derive $N_{L/M_{N}^{k}}\simeq O_{L}(-1)^{\oplus k-N+2}\oplus O_{L}^{\oplus
{2N-k-5}}$. 
In Section 2, we study
the moduli space $\overline{M}_{0,0}(\bP^{1},2)$ from the point of view of
stability and identify it with $\bP^{2}$ and show that
the moduli space $\overline{M}_{0,0}(\bP^{1},2)$ is isomorphic to
$\bP(S^{2}Q)$ over $G$.
In section 4, we describe $H^{1}(\phi^{*}N_{L/M_{N}^{k}})$ as an coherent
sheaf over $\bP(S^{2}Q)^{stack}$.
In section 5, we
derive the main theorem (\ref{re}) of this paper by using Segre classes. In
Section 6, we mention
some generalization to  degree 3 Gromov-Witten invariants.

\section{Lines on a hypersurface}
\label{sec:lines on a hypersurface}
Let $M$ be a generic hypersurface of degree $k$ of the projective space $\bP^{N-1}=\bP(V)$. 
We assume $2N-5\geq k\geq N-2\geq 2$ throughout this note. 
In this note we count the number of rational curves of virtual degree two, namely rational curves 
which doubly cover lines on $M$. 
\par
Let $\bP=\bP(V)$ be the projective space 
parameterizing all one-dimensional quotients of $V$, 
which is usually denoted by $\bP(V)$ in 
the standard notation in algebraic geometry. 
In this notation let $W$ be a subspace of $V$. Then $\bP(W)$ is 
naturally a linear subspace of $\bP(V)$ 
of dimension $\dim W - 1$.\par
Let $G(2,V)$ be the Grassmann variety of lines in $\bP(V)$, 
the scheme parameterizing all lines of $\bP=\bP(V)$. 
This is also the universal scheme parameterizing all 
one-dimensional quotient linear spaces of $V$. 
Let $W$ be a two dimensional quotient linear space, 
$\psi\in G(2,V)$, namely 
$\psi : \bP(W)\to \bP(V)$ the natural immersion 
and $i^*_{\psi}:V\to W$ the quotient homomorphism. 
The space $W$ is denoted by $W(\psi)$ when necessary.\par
There exists the universal bundle $Q_{G(2,V)}$ over $G(2,V)$ and a 
homomorphism $i^{\op{univ}*}: O_{G(2,V)}\otimes V\to Q_{G(2,V)}$
whose fiber $i^{\op{univ}*}_{\psi}:V\to Q_{G(2,V),\psi}$ is the quotient 
$i^*_{\psi}:V\to W(\psi)$ of $V$ corresponding to $\psi$.

\subsection{Existence of a line on $M$}
Let $L=\bP(W)$ be a line of $\bP$, equivalently $W\in G(2,V)$. Then the condition 
$L\subset M$ imposes at most $k+1$ conditions on $W$, 
while the number of moduli of lines of $\bP$ equals 
$\dim G(2,V)=2N-4$.  Hence we infer 
\begin{lemma}\label{lemma:existence of a line on M}
If $2N\geq k+5$, then there exists at least a line on $M$.
\end{lemma}

See also [Katz,p.152]. Let $G$ be the subscheme of $G(2,V)$ 
parameterizing all lines of $\bP(V)$ lying on $M$, $Q=(Q_{G(2,V)})_{|G}$ 
the restriction of $Q_{G(2,V)}$ to $G$. By Lemma~\ref{lemma:existence of a line on M}, 
$G$ is nonempty.
Let $i^*:O_{G}\otimes V\to Q$ be the restriction of $i^{\op{univ}*}$ to $G$. 
Let $P=\bP(Q)$ and $\pi : P\to G$ the natural projection. Then $\pi$
is the universal line of $M$ over $G$, to be more exact, 
the universal family  over $G$ of lines lying on $M$. In other words,  
the natural epimorphism 
$i^* : O_{G}\otimes V\to Q$
induces a morphism $i:P\to \bP_G(V):=G\times \bP(V)$, which is 
a closed immersion into $\bP_G(V)$, thus $P$ is a subscheme of $\bP_G(V)$ 
such that $\pi=(p_1)_{|P}$.  Let $L_{\psi}=\bP(Q_{\psi})$.  Note that
\begin{equation*}
L_{\psi}=P_{\psi} :=\pi^{-1}(\psi)\simeq \bP(Q_{\psi})
\subset \{\psi\}\times \bP(V)\simeq \bP(V).
\end{equation*}

\subsection{The normal bundle $N_{L/M}$}
The argument of this section is standard and well known.
Let $\bP=\bP(V)$, $L=\bP(W)$ and $i^*_W:V\twoheadrightarrow W\in G$. 
 Let us recall the following exact sequence:
\begin{equation*}
\CD
0 @>>> O_{\bP} @>>>  O_{\bP}(1)\otimes V^{\vee}@>D>>
T_{\bP}@>>>  0
\endCD
\end{equation*}
where 
the homomorphism $D$ is defined by 
\begin{align*}
D(a\otimes v^{\vee}):&=aD_{(v^{\vee})}\quad (a\in O_{\bP}(1))\\
(D_{v^{\vee}}F)(u^{\vee}):&=(\frac{d}{dt} F(u^{\vee}+tv^{\vee}))_{t=0}
\end{align*}
for a homogeneous polynomial 
$F\in S(V)$ and $u^{\vee}, v^{\vee}\in V^{\vee}$. 
We note
$H^0(O_{\bP}(1))\otimes V^{\vee}=V\otimes V^{\vee}=\End(V,V)$
and that the image of $H^0(O_{\bP})$ in 
$\End(V,V)$ is $\bC\op{id_V}$.
We also have the following exact sequences:
\begin{equation*}
\CD
0 @>>> T_{L} @>>>
(T_{\bP})_L@>>> N_{L/\bP} @>>> 0\\
0 @>>> O_L @>>>  O_L(1)\otimes V^{\vee}@>D_L>>
(T_{\bP})_L@>>>  0.
\endCD
\end{equation*}

\begin{lemma}\label{lemma:normal bundle of L in P}Let $L=\bP(W)$. Then 
\begin{equation*}
N_{L/\bP}\simeq O_L(1)\otimes (V^{\vee}/W^{\vee}),\ 
H^0(N_{L/\bP})\simeq W\otimes (V^{\vee}/W^{\vee}).
\end{equation*}
\end{lemma}
\begin{pf}
The assertion is clear from the following commutative diagram with exact rows and columns:
\begin{equation*}
\CD
0 @>>> O_L @>>>  O_L(1)\otimes W^{\vee}@>(D_L)_{|W^{\vee}}>>
T_L @>>>  0\\
@VV{}V @VV{}V @VV{\op{id}\otimes i^{\vee}}V @VV{}V@VVV \\
0 @>>> O_L @>>>  O_L(1)\otimes V^{\vee}@>D_L>>
(T_{\bP})_L@>>>  0\\
@VVV @VV{}V @VV{}V @VVV @VVV \\
0 @>>> 0 @>>>
O_L(1)\otimes (V^{\vee}/W^{\vee})@>>> N_{L/\bP} @>>> 0
\endCD
\end{equation*}

The second assertion is clear from $H^0(L,O_L(1))=W$. \qed
\end{pf}

Since $T_L\simeq O_L(2)$, there follow exact sequences
\begin{equation*}
\CD
0 @>>> H^0(T_{L}) @>>>
H^0((T_{\bP})_L)@>>> H^0(N_{L/\bP}) @>>> 0\\
0 @>>> H^0(O_L) @>>>  H^0(O_L(1))\otimes V^{\vee}@>H^0(D_L)>>
H^0((T_{\bP})_L)@>>>  0.
\endCD
\end{equation*}

  We also note 
\begin{equation*}
H^0(T_{L})=\op{Lie}\op{Aut}^0(L)
=\op{End}(W,W)/\op{center}=\op{End}(W,W)/\bC\op{id_W}.
\end{equation*}

Since $H^0(O_L(1))=W$, we see 
\begin{equation*}
H^0((T_{\bP})_L)=W\otimes V^{\vee}/\op{Im}H^0(O_L)=\op{Hom}(V,W)/\bC i^*_W.
\end{equation*}

Hence we again see
\begin{align*}
H^0((N_{L/\bP}))&=(\op{Hom}(V,W)/\bC i^*_W)/(\op{Hom}(W,W)/\bC \op{id}_W)\\
&=W\otimes (V^{\vee}/W^{\vee})=\op{Hom}(V/W,W).
\end{align*}

For any line $L=\bP(W)$ of $\bP$ 
the following sequence is exact:
\begin{equation}\label{eq:normal sequence}
0\to N_{L/M}\to N_{L/\bP} \to (N_{M/\bP})_{L}(\simeq O_L(k)) \to 0.
\end{equation}

Hence so is the following sequence as well:
\begin{equation*}
\CD
0@>>>H^0(N_{L/M})@>>>H^0(N_{L/\bP})@>H^0(D_L)>>H^0(O_L(k))\\
@>>> H^1(N_{L/M})@>>> 0.
\endCD
\end{equation*}

Hence we have
\begin{lemma}\label{lemma:normal sequence}¡¡The following is exact:
\begin{equation}
\CD
0\to H^0(N_{L/M})\to W\otimes (V^{\vee}/W^{\vee})
\overset{H^0(D_L)}{\longrightarrow} S^kW\to H^1(N_{L/M})\to  0.
\endCD
\end{equation}
\end{lemma}

\begin{cor}$\dim G\geq 2N-k-5$, equality holding if $H^1(N_{L/M})=0$.
\end{cor}
\begin{pf}As is well-known, $\dim G\geq h^0(N_{L/M})-h^1(N_{L/M})$. 
Note $\dim W\otimes (V^{\vee}/W^{\vee})=2(N-2)$ and $\dim S^kW=k+1$. Hence the corollary follows from 
Lemma~\ref{lemma:normal sequence}. \qed
\end{pf}

\begin{lemma}\label{lemma:normal bundle formula}
For a generic line $L$ on a generic hypersurface $M$ of degree $k$ 
\begin{enumerate}
\item[(i)]\ $N_{L/M}\simeq O_L^{\oplus a}\oplus O_L(-1)^{\oplus b}$, 
where $a=2N-k-5$ and $b=k-N+2$,
\item[(ii)]\ $\coker H^0(D^{-}_L)\simeq S^{k-1}W/(V^{\vee}/W^{\vee})$
where $D^{-}_L:=D_L\otimes O_L(-1)$.
\end{enumerate}
\end{lemma}
\begin{pf}Let $M$ be a generic hypersurface of degree $k$ and $L$ a generic line $L$ on $M$. 
Without loss of generality we may assume that 
$W^{\vee}$ is generated by $e_1^{\vee}$ and $e_2^{\vee}$, 
in other words, $\psi : L\to \bP$ is given by
$$\psi : [s:t]\to [x_1,\cdots,x_N]=[s,t,0,\cdots,0].
$$

Then $F$, the polynomial of degre $k$ defining $M$, is written as 
$$F=x_3F_3+x_4F_4+\cdots+x_NF_N
$$
for some polynomials $F_j$ of degree $k-1$. Let $f_j=\psi^*F_j=F_j(s,t,0,\cdots,0)$.
\par
Now we consider the exact sequence 
\begin{equation*}
\CD
0@>>>H^0(N_{L/M}(-1))@>>>H^0(N_{L/\bP}(-1))@>H^0(D^{-}_L)>>H^0(O_L(k-1))\\
@>>> H^1(N_{L/M}(-1))@>>> 0.
\endCD
\end{equation*}
where we note
$H^0(N_{L/\bP}(-1))=V^{\vee}/W^{\vee}$.  Hence the following is exact:
\begin{equation}
\CD
0@>>>H^0(N_{L/M}(-1))@>>>V^{\vee}/W^{\vee}@>H^0(D^{-}_L)>>S^{k-1}W\\
@>>> H^1(N_{L/M}(-1))@>>> 0
\endCD
\end{equation}
where $H^0(D^{-}_L)$ is given by
$H^0(D^{-}_L)(e_j^{\vee})=f_j\quad (j=3,4,\cdots,N)$.\par

A generic choice of $F$ implies a generic choice of 
degree $k-1$ polynomials $f_j$ $(j=3,4,\cdots,N)$ in $s$ and $t$. By the assumptions  
\begin{gather*}
\dim S^{k-1}W=k\geq N-2=\dim V^{\vee}/W^{\vee},\\  
\dim W\otimes V^{\vee}/W^{\vee}=2(N-2)\geq k+1=\dim S^{k}W,
\end{gather*}
the generic choice of $F$ implies that 
we can choose  $f_j\in S^{k-1}W$ $(j=3,4,\cdots,N)$ (and fix once for all) such that
\begin{enumerate}
\item[(iii)] $f_j$ $(j=3,4,\cdots,N)$ are linearly independent,
\item[(iv)] $Wf_3+Wf_4+\cdots+Wf_N=S^kW$.
\end{enumerate}

Hence $H^0(D^{-}_L)$ is injective by (iii). 
It follows that $H^0(N_{L/M}(-1))=0$. Hence (ii) is clear.  
Next we consider $H^0(D_L)$.  By (iv), we see 
$$S^kW=W\cdot H^0(D^{-}_L)(V^{\vee}/W^{\vee})=H^0(D_L)(W\otimes V^{\vee}/W^{\vee}),$$
whence $H^0(D_L)$ is surjective. It follows that $H^1(N_{L/M})=0$. 
Hence $N_{L/M}\simeq O_L^{\oplus a}\oplus O_L(-1)^{\oplus b}$ for some $a$ and $b$.
Since $a+b=\rank(N_{L/M})=N-3$ and $-b=\deg(N_{L/M})=N-2-k$, we have (i).\qed
\end{pf}

\subsection{Lines on a quintic hypersurface in $\bP^4$}
See [Katz, Appendix A]  for the subsequent examples. 
Let $N=5$ and $k=5$. Hence $M$ is a 
hypersurface of degree 5 in $\bP^4$, a Calabi-Yau 3-fold. 
Let  
$$F = x_4x_1^4 + x_5x_2^4 + x_3^5 + x_4^5 + x_5^5.
$$

First we note that $M=\{F=0\}$ is nonsingular. Let $L=\{x_3=x_4=x_5=0\}=\{[s,t,0,0,0]\}$. 
In this case $f_3=0$, $f_4=s^4$  and $f_5=t^4$. In the exact sequence (1) we see
$H^0(N_{L/M}(-1))=\Ker H^0(D^{-}_L)=\bC e_3^{\vee}$ and 
$H^1(N_{L/M}(-1))=\coker H^0(D^{-}_L)$ is 3-dimensional. 
Hence  $N_{L/M}=O_L(1)\oplus O_L(-3)$.
\par
We summarize the above. If $\dim\Ker H^0(D^{-}_L)=1$ and if $M$ is nonsingular,
then $N_{L/M}=O_L(1)\oplus O_L(-3)$. 
Hence $H^0(N_{L/M})=\Ker H^0(D_L)=W\otimes \Ker H^0(D^{-}_L)$ is 2-dimensional. 
Therefore we can choose $f_3=0$ and a linearly independent pair $f_4$ and $f_5\in S^4W$ so that
$Wf_4+Wf_5$ is 4-dimensional. The choice $f_4=s^4$ and $f_5=t^4$ satisfies the conditions.
This enables us to find a nonsingular hypersurface $M$ as above. 
However if we choose $f_3=0$, $f_4=s^4$ and $f_5=s^3t$, then $Wf_4+Wf_5$ is 3-dimensional. 
Hence $M$ is singular.\par 
Next in the same manner we find $L$ on a nonsingular hypersurface 
$M$ with $N_{L/M}=O_L\oplus O_L(-2)$
or $N_{L/M}=O_L(-1)^{\oplus 2}$. Let 
$$F = x_3x_1^4+x_4x_1^3x_2+x_5x_2^4+x_3^5+x_4^5+x_5^5.
$$

Then we have $f_3=s^4$, $f_4=s^3t$ and $f_5=t^4$. Since
$Wf_3+Wf_4+Wf_5$ is 5-dimensional, 
$H^0(N_{L/M}(-1))=\Ker H^0(D^{-}_L)=0$, $H^0(N_{L/M})=\Ker H^0(D_L)=\bC(te_3^{\vee}-se_4^{\vee})$. 
We see also 
that $\dim H^1(N_{L/M})=\dim\coker H^0(D_L)=1$ and $N_{L/M}=O_L\oplus O_L(-2)$. 
The hypersurface 
$M=\{F=0\}$ is easily shown to be nonsingular. \par
If $F = x_3x_1^4+x_4x_1^2x_2^2+x_5x_2^4+x_3^5+x_4^5+x_5^5$ and $M=\{F=0\}$, 
then $N_{L/M}=O_L(-1)^{\oplus 2}$.

\subsection{Lines on a generic hypersurface $M_{7}^{8}$ of $\bP^6$}
\label{subsec:Lines on a generic hypersurface}
Let $N=7$ and $k=8$. In view of Lemma~\ref{lemma:existence of a line on M} 
there exists a line $L$ on any generic hypersurface 
$M$ of degree 8 in $\bP(V)=\bP^6$. In view of Lemma~\ref{lemma:normal bundle formula},
$a=1$, $b=3$ and $N_{L/M}\simeq O_L\oplus O_L(-1)^{\oplus 3}$.
For example let $L: x_j=0\ (j\geq 3)$ and we take 
\begin{gather*}
F_3=8x_1^7, F_4=8x_1^6x_2, F_5=8x_1^4x_2^3, F_6=8x_1^2x_2^5, F_7=8x_2^7,\\ 
F=x_3F_3+x_4F_4+x_5F_5+x_6F_6+x_7F_7+x_3^8+x_4^8+x_5^8+x_6^8+x_7^8.
\end{gather*}
and let $M=M_{7}^{8}:F=0$.
We see that $M$ is nonsingular near $L$ and 
has at most isolated singularities. However 
it is still unclear to us 
whether $M=M_{7}^{8}$ is nonsingular everywhere. 
The space $H^0(N_{L/M})$ is spanned by $te_3^{\vee}-se_4^{\vee}$, hence 
an infinitesimal deformation $L_{\varep}$ of $L$ is given by 
\begin{equation*}[s,t]\mapsto [s,t,\varep t, -\varep  s,0,0,0]
\end{equation*} 
which yields $F_{|L_{\varep}}=\varep^8(s^8+t^8)\equiv 0\ \mod \varep^8$. 
Since $H^1(N_{L/M})=0$, this infinitesimal deformation is integrable and  
$G$ ($:=$ the moduli of lines of $\bP^6$ contained in $M$) 
is nonsingular and one dimensional at the point $[L]$. \par
We note that $M$ also contains 8 lines
$$L':=L'_{\varep_8} : \varep_8x_1-x_2=x_3+\varep_8x_4=x_j=0\ (j\geq 5),
$$ 
with $N_{L'/M}=O_{L'}(1)^{\oplus 3}\oplus O_{L'}(-6)$ 
where $\varep_8^8=-1$.

\section{Stability}\label{sec:stability}

\begin{defn}Suppose that a reductive algebraic group 
$G$ acts on a vector space $V$.
 Let $v\in V$, $v\neq 0$. 
\begin{enumerate}
\item[(1)] the vector $v$ is said to be {\it semistable} if there exists a
$G$-invariant homogeneous polynomial $F$ on $V$
such that $F(v)\neq 0$,
\item[(2)] the vector $v$ is said to be {\it stable} 
if $p$ has a closed $G$-orbit in $X_{ss}$
and the stabilizer subgroup
of $v$ in $G$ is finite.
\end{enumerate}
Let $\pi:V\setminus\{0\}\to
\bP(V^{\vee})$ be the natural surjection.
Then $v\in V$ is semistable (resp. stable) if and only if
$\pi(v)$ is semistable (resp. stable).
\end{defn}

\subsection{Grassmann variety}\label{subsec:Grassman}
Let $V$ be an $N$-dimensional vector space, and $G(r,N)$ the Grassmann variety 
parameterizing all $r$-dimensional quotient spaces of $V$. 
Here is a natural way of understanding 
$G(r,N)$ via GIT-stability. Let $U$ be an $r$-dimensional vector space, 
$X=\Hom(V,U)$ and  $\pi:X\setminus\{0\}\to \bP(X^{\vee})$ the natural map. 
Then $\SL(U)$ acts on $X$ from the left by:
\begin{equation*}
(g\cdot\phi^*)(v)=g\cdot(\phi^*(v))\quad \op{for}\ \phi^*\in X,\ v\in V. 
\end{equation*}
We see that for $\phi^*\in X$
\begin{center}
 $\phi^*$ is $\SL(U)$-stable $\Longleftrightarrow$ $\rank\phi^*=r$,\\
 $\phi^*$ is $\SL(U)$-semistable
 $\Longleftrightarrow$ $\phi^*$ is $\SL(U)$-stable.
\end{center}

In fact, if $\rank\phi^*= r-1$, then there is a one-parameter torus $T$ of 
$\SL(U)$ such that the closure of the orbit 
$T\cdot\phi$ contains the zero vector as 
the following simple example $(r=2)$ shows
\begin{equation*}
\lim_{t\to 0}\begin{pmatrix}t&0\\
0&t^{-1}
\end{pmatrix}
\begin{pmatrix}
a_{11}&a_{12}&\cdots&a_{1N}\\
0&0&\cdots&0
\end{pmatrix}
=\lim_{t\to 0}\begin{pmatrix}
ta_{11}&ta_{12}&\cdots&ta_{1N}\\
0&0&\cdots&0
\end{pmatrix}.
\end{equation*}

Let $X_{\op{s}}$ be the set of all (semi)stable points and 
$\bP_s$ the image of $X_{\op{s}}$ by $\pi$. 
It is, as we saw above, just the set of all $\phi\in X$ with $\rank\phi^*=r$.
Therefore 
the GIT-orbit space $\bP_s\quot\SL(U)$ is 
the orbit space $\bP_s/\SL(U)$ by the free action, 
the Grassmann variety $G(r,N)$.

\subsection{Moduli of double coverings of $\bP^1$ (1)}
\label{subsec:moduli of double coverings 1}
Let $W$ and $U$ be a pair of two dimensional vector spaces, $X=\Hom(W,S^2U)$, and 
$\pi:X\setminus \{0\}\to\bP(X^{\vee})$ the natural morphism.   
Note that $\SL(U)$ acts on $S^2U$ from the left via the natural action:
$\sigma(u_1u_2)=\sigma(u_1)\sigma(u_2)$ for $\forall u_1,u_2\in U$.
Thus $\SL(U)$ acts on $X$ from the left 
in the same manner in the subsection~\ref{subsec:Grassman}. 
\begin{lemma}Let $\phi^*\in X$.\par
\begin{enumerate}
\item[(i)]
 $\phi^*$ is unstable iff $\phi^*(w)$ has a double root for any $w\in W$,
\item[(ii)] $\phi^*$ is semistable  iff  $\phi^*(w)$ has no double roots 
 for some nonzero $w\in W$,
\item [(iii)] $\phi^*$ is stable
  iff  $\phi^*(W)$ is a base-point free linear subsystem of $S^2U$ on $\bP(U)$.
\end{enumerate}
\end{lemma}
\begin{proof}We note that $\phi^*$ is unstable iff there is a suitable basis $s$ and $t$ of $U$ such that 
$\phi^*(w)=a(w)s^2$ for any $w\in W$ since a torus orbit $T\cdot \phi^*$ contains the zero vector. 
This proves (i). This also proves (ii). Next we prove (iii). If $\phi^*(W)$ has a base point, then 
it is clear that $\phi^*$ is not stable. If $\phi^*$ is semistable and it is not stable, then 
we choose a basis $s$, $t$ of $U$ and a basis $w_1$, $w_2$ of $W$ such that $\phi^*(w_1)=st$. If 
$\phi(w_2)=as^2+bst$, then $\phi^*$ is not stable. This proves the lemma. 
\end{proof}

\begin{thm}\label{thm:quotient}
Let $X_{ss}$ be the Zariski open subset of $X$ consisting of all semistable points of $X$, 
$\pi(X_{ss})$ the image of $X_{ss}$ by $\pi$,
and $Y:=\pi(X_{ss})\quot\SL(U)$.  Then $Y\simeq\bP^2$. 
\end{thm}
\begin{pf}First consider a simplest case.  We choose a basis $s$, $t$ of $U$. 
Let $w_1$ and $w_2$ be a basis of $W$, $T$ the subgroup of $\SL(U)$ 
of diagonal matrices and $X'=\{\phi^*\in X; \phi^*(w_1)=2st\}$. 
Let $Z'=\SL(U)\cdot X'$. \par 
We note that $Z'$ is an $\SL(U)$-invariant subset of $X_{ss}$. 
We prove $\pi(Z')\quot\SL(U)\simeq\bC^2$. 
Let $\phi^*$ and $\psi^*$ be points of $X'$. 
Let $\phi^*(w_2)=As^2+2Bst+Ct^2$ and $\psi^*(w_2)=as^2+2bst+ct^2$. 
Then it is easy to check 
\begin{align*}
g\cdot\phi^*=\psi^*\ \op{for}\ \exists\ g\in \SL(U)
&\Longleftrightarrow   
 g\cdot\phi^*=\psi^*\ \op{for}\ \exists\ g\in T\\
&\Longleftrightarrow  A=au^2,\ B=b,\ C=u^{-2}c\ \op{for}\ \exists\ u\neq 0.
\end{align*}
Therefore each equivalence class of $\pi(Z)\quot\SL(U)$ is represented by the pair 
$(AC,B)$, which proves $\pi(Z)\quot\SL(U)\simeq\bC^2$.\par
Now we prove the lemma.  
Let $\phi^*\in X_{ss}$,
$\phi_j=\phi^*(w_j)$ and $\phi_0=-(\phi_1+\phi_2)$. 
Let
\begin{align*}
\phi_0&=r_1s^2+2r_2st+r_3t^2,\\
\phi_1&=p_1s^2+2p_2st+p_3t^2,\\
\phi_2&=q_1s^2+2q_2st+q_3t^2,
\end{align*}
and we define
\begin{gather*}
D_1=p_2^2-p_1p_3,\quad D_2=q_2^2-q_1q_3,\\
D_0=r_2^2-r_1r_3=D_1+D_2+2p_2q_2-(p_1q_3+p_3q_1).
\end{gather*}

To show the lemma, we prove the more precise isomorphism
\begin{equation*}\pi(X_{ss})\quot\SL(U)=\Proj\bC[D_0,D_1,D_2]
\end{equation*} 
For this purpose we define 
$Y_j=\pi(\{\phi^*\in X_{ss};\phi_j\ \op{has\ no\ double\ roots}\})\quot\SL(U)$.
It suffices to prove $Y_1=
\Spec\bC[\frac{D_0}{D_1},\frac{D_2}{D_1}]$ by reducing it to the first simplest case.\par

Let $\phi^*\in Y_1$. Let $\alpha$ and $\beta$ be the roots of $\phi_1=0$. 
By the assumption $\phi_1$ has no double roots, hence $\alpha\neq\beta$. 
Let 
\begin{gather*}
u=\frac{1}{\gamma}(s - \alpha t),\quad v = \frac{1}{\gamma}(s - \beta t),\quad
g=\frac{1}{\gamma}\begin{pmatrix}
1  & - \alpha\\
1  & - \beta
\end{pmatrix}
\end{gather*}
where $\gamma=\sqrt{\alpha-\beta}$. Note that $g\in\SL(U)$. 
Hence we see
\begin{align*}
(\phi_1(s,t),\phi_2(s,t)) &\equiv (p_1\gamma^4 uv, A_1u^2+2B_1uv+C_1v^2)
\end{align*}
where 
\begin{align*}
A_1&=q_1^2\beta^2+2q_2\beta+q_3,\\
-B_1&=q_1\alpha\beta+q_2(\alpha+\beta)+q_3,\\
C_1&=q_1^2\alpha^2+2q_2\alpha+q_3.
\end{align*}

Thus we see
\begin{align*}
(\phi_1(s,t),\phi_2(s,t)) &\equiv (2st, As^2+2Bst+Ct^2) 
\end{align*}
where
\begin{gather*}
A=\frac{2A_1}{p_1\gamma^4},\
B=\frac{2B_1}{p_1\gamma^4},\
C=\frac{2C_1}{p_1\gamma^4},\
p^2_1\gamma^4=4D_1,\\
AC=B^2-\frac{D_2}{D_1},\ B=\frac{D_0-D_1-D_2}{2D_1}.
\end{gather*}

Therefore by the first half of the proof
\begin{equation*}
Y_1\simeq \Spec\bC[AC,B]=\Spec\bC[\frac{D_0}{D_1},\frac{D_2}{D_1}].
\end{equation*}
\par
This completes the proof of the lemma.\qed
\end{pf}

\begin{cor}\label{cor:conic}
Let $Y^s=\pi(X_s)\quot\SL(U)$. 
Then $Y\setminus Y^s$ is a conic of $Y$ defined by
\begin{equation*}Y\setminus Y^s:D_0^2+D_1^2+D_2^2-2D_0D_1-2D_1D_2-2D_2D_0=0.
\end{equation*}
\end{cor}
\begin{pf}In view of Theorem~\ref{thm:quotient}, $Y_1\simeq\Spec\bC[AC, B]$. 
The complement of $Y_s$ in $Y_1$ is then the curve defined by $AC=0$, 
which is easily identified with the above conic.\qed
\end{pf}

\begin{cor}\label{cor:quotient of X^0}
Let $X^0$ be the Zariski open subset of $X$ 
consisting of all semistable points $\phi^*$ of $X$ with $\rank\phi^*=2$, 
and let $Y^0:=\pi(X^0)\quot\SL(U)$.  Then  
 $Y^0\simeq \pi(X^0)/\SL(U)\simeq Y\simeq\bP^2$. 
\end{cor}
\begin{pf}It suffices to compare $Y_1$ and $Y^0\cap Y_1$. As in the proof of 
Theorem~\ref{thm:quotient} we let $X'=\{\phi^*\in X; \phi^*(w_1)=2st\}$. 
Let $Z=\SL(U)\cdot X'$ and $Z^0=\SL(U)\cdot (X'\cap X^0)$. \par 
Then with the notation in Theorem~\ref{thm:quotient}, we recall
$X'=\{\phi^*\in X; \phi^*(w_1)=2st, \phi^*(w_2)=As^2+2Bst+Ct^2\}$, 
$\pi(Z)\quot\SL(U)\simeq\Spec\bC[AC,B]$ where
\begin{equation*}
X'\cap X^0
=\{\phi^*\in X'; A\neq 0\ \op{or}\ C\neq 0\}.
\end{equation*}

In the same manner as before we see $\pi(Z^0)\quot\SL(U)\simeq\Spec\bC[AC,B]$, 
whence $\pi(Z^0)\quot\SL(U)=\pi(Z)\quot\SL(U)$. This proves $Y^0\cap Y_1=Y_1$.
This completes the proof of the corollary.\qed
\end{pf}

\subsection{Moduli of double coverings of $\bP(W)$\ (2)}
There is an alternative way of understanding $\pi(X_{ss})\quot\SL(U)\simeq \bP^2$
by using the isomorphism $S^2\bP^1\simeq\bP^2$. We use the following convention to denote 
a point of $\bP(U)=U^{\vee}\setminus\{0\}/\bG_m$:
$(u:v)=us^{\vee}+vt^{\vee}\in U^{\vee}$ 
where $s^{\vee}$ and $t^{\vee}$ are a basis dual to $s$ and $t$.
In what follows we fix a basis $w_1$ and $w_2$ of $W$. 
Let $P:=(a_1:a_2)$ and $Q:=(b_1:b_2)$ be a pair of 
points of $\bP(W)\simeq \bP^1$. If $P\neq Q$, 
there is a double covering $\phi:\bP(U)\to\bP(W)$ 
ramifying at $P$ and $Q$, unique up to isomorphism once we fix the base $w_1$ and $w_2$: 
\begin{equation*}
\frac{b_2w_1-b_1w_2}{a_2w_1-a_1w_2}
=(\frac{t}{s})^2.
\end{equation*}
Thus $\phi$ is given explicitly by 
\begin{equation*}\phi_1:=\phi^*(w_1)=b_1s^2-a_1t^2,\ 
\phi_2:=\phi^*(w_2)=b_2s^2-a_2t^2,\ \phi_0=-(\phi_1+\phi_2)
\end{equation*}
for which we have  
\begin{equation*}D_1=a_1b_1,\ D_2=a_2b_2,\ D_0=(a_1+a_2)(b_1+b_2).
\end{equation*}

The isomorphism $S^2\bP^1\simeq\bP^2$ is given by $(P,Q)\mapsto (D_0,D_1,D_2)$. 
This shows 
\begin{cor} We have a natural isomorphism:
$Y\simeq\bP(S^2W)$. 
\end{cor}

\section{The virtual normal bundle of a double covering}
\label{sec:virtual normal bundle}
\subsection{The case $N=7$ and $k=8$ revisited} We revisit  
the example 
in the subsection~\ref{subsec:Lines on a generic hypersurface}. 
Let $N=7$ and $k=8$. Let $L: x_j=0\ (j\geq 3)$ and we take 
\begin{gather*}
F_3=8x_1^7, F_4=8x_1^6x_2, F_5=8x_1^4x_2^3, F_6=8x_1^2x_2^5, F_7=8x_2^7,\\ 
F=x_3F_3+x_4F_4+x_5F_5+x_6F_6+x_7F_7+x_3^8+x_4^8+x_5^8+x_6^8+x_7^8.
\end{gather*}
and let $M=M^5_8:F=0$. We often denote $L$ also by $\bP(W)$ 
with $W$ a two dimensional vector space for later convenience.
Since $H^0(D^{-}_L)$ is injective and $H^0(D_L)$ is surjective, 
we have $N_{L/M}\simeq O_L\oplus O_L(-1)^{\oplus 3}$. Hence 
$H^1(N_{L/M}(-1))=H^1(O_L(-2)^{\oplus 3})$ is 3-dimensional. 
As we see easily, this follows also 
from the fact that $\coker H^0(D^{-}_L)$ is freely generated by 
$x_1^5x_2^2$, $x_1^3x_2^4$ and $x_1x_2^6$. \par
Let $\phi^*=(\phi_1,\phi_2)\in X^0$. Then
$\Ker H^0(\phi^*D_L)$ is generated 
by a single element $\phi_2e_3^{\vee}-\phi_1e_4^{\vee}$,
while $\coker H^0(\phi^*D_L)$ is generated by 
$S^2U\cdot\phi_1^5\phi_2^2$, $S^2U\cdot\phi_1^3\phi_2^4$ 
and $S^2U\cdot\phi_1\phi_2^6$.
To be more precise, we see 
\begin{equation*}
\coker H^0(\phi^*D_L)=\{\phi_1^5\phi_2^2,\phi_1^3\phi_2^4,\phi_1\phi_2^6\}\otimes
S^2U/\{\phi_1,\phi_2\}.
\end{equation*}
In fact, this is proved as follows: first we consider the case where 
$\phi_1$ and $\phi_2$ has no common zeroes. In this case $\phi^*$ gives rise to a 
double covering $\phi:\bP(U)\to \bP(W)\ (=L)$, which we denote by 
$L_{\phi}$ for brevity. By pulling back by $\phi^*$
the normal sequence $0\to N_{L/M}\to N_{L/\bP} \to O_L(k) \to 0$ $(k=8)$ for the line $L$
we infer an exact sequence 
\begin{equation*}
0\to \phi^*N_{L/M}\to \phi^*N_{L/\bP} \overset{\phi^*D_L}\longrightarrow \phi^*O_L(k) \to 0, 
\end{equation*}
which yields an exact sequence 
\vskip -0.5cm
\begin{equation*}
\CD
0 @>>> H^0(\phi^*N_{L/M}) @>>>  S^2U\otimes (V^{\vee}/W^{\vee})@>H^0(\phi^*D_L)>>
H^0(O_{L_{\phi}}(2k))  \\
@>>> H^1(\phi^*N_{L/M}) @>>> 0. 
\endCD
\end{equation*}

Let 
$\eta=q_3e_3^{\vee}+\cdots+q_7e_7^{\vee}\in\Ker H^0(\phi^*D_L)$, $q_j\in S^2U$. 
Then we have
\begin{equation*}\phi_1^2(q_3\phi_1^5+q_4\phi_1^4\phi_2+q_5\phi_1^2\phi_2^3+q_6\phi_2^5)=-q_7\phi_2^7.
\end{equation*}

Since $\phi_1$ and $\phi_2$ are mutually prime and $q_j$ is of degree two, we have $q_7=0$ and 
\begin{equation*}\phi_1^2(q_3\phi_1^3+q_4\phi_1^2\phi_2+q_5\phi_2^3)=-q_6\phi_2^5,
\end{equation*}

Hence $q_6=0$ and similarly we infer also $q_5=0$. Thus we have $q_3\phi_1+q_4\phi_2=0$. 
This proves that $\Ker H^0(\phi^*D_L)$ is generated by $\phi_2e_3^{\vee}-\phi_1e_4^{\vee}$.\par
Next we prove that $\coker H^0(\phi^*D_L)$ is generated by 
$\phi^*\coker H^0(D^{-}_L)$ over $S^2U$, in fact over $S^2U/\phi^*(W)$.
Without loss of generality we may assume that $\phi_1=2st$ and 
$\phi_2=\lambda s^2+2\nu st + t^2$ for some $\lambda\neq 0$ and $\nu\in\bC$. 
Let $\phi^*W=\{\phi_1,\phi_2\}$.
Then one checks $U\cdot\phi^*W=S^3U$, and hence $S^2U\cdot\phi^*W=S^4U$,
$S^{2m-2}U\cdot\phi^*W=S^{2m}U$ for $m\geq 2$. It follows $S^2U\cdot\phi^*(S^{m-1}W)=S^{2m}U$ 
for $m\geq 1$. In fact, by the induction on $m$
\begin{align*}S^2U\cdot\phi^*(S^mW)&=S^2U\cdot \phi^*(W)\cdot\phi^*(S^{m-1}W)\\
&=S^4U\cdot\phi^*(S^{m-1}W)\\
&=S^2U\cdot (S^2U\cdot\phi^*(S^{m-1}W))\\
&=S^2U\cdot S^{2m}U=S^{2m+2}U.
\end{align*}

Therefore $H^0(O_{L_{\phi}}(2k))=S^{16}U=S^2U\cdot\phi^*(S^7W)$.  Hence 
\begin{align*}
\coker H^0(\phi^*D_L)&=S^{16}U/\image H^0(\phi^*D_L)\\
&=S^2U\cdot\phi^*(S^7W)/S^2U\cdot\phi^*(\image H^0(D^{-}_L))\\
&=(S^2U/\phi^*(W))\cdot\phi^*(S^7W/\image H^0(D^{-}_L)).
\end{align*}
because $\coker H^0(D^{-}_L)=S^7W/\image H^0(D^{-}_L)$ and 
$W\cdot S^7W\subset W\cdot \image H^0(D^{-}_L))=S^8W$ by the choice of $L$.
This proves that $\coker H^0(\phi^*D_L)$ is generated 
by $\phi^*\coker H^0(D^{-}_L)$ over $S^2U/\phi^*(W)$. It follows 
$\coker H^0(\phi^*D_L)
=(\phi^*\coker H^0(D^{-}_L))\otimes (S^2U/\phi^*W)$.\par
Finally we consider the case where $\phi_1$ and $\phi_2$ has a common zero. 
In this case we may assume $\phi_1=2st$ and $\phi_2=2\nu st+t^2$. In this case
$L_{\phi}$ is a chain of two rational curves $C'_{\phi}$ and $C''_{\phi}$ where 
$C_{\phi}$ is the proper transform of $\bP(U)$, where 
the double covering map from $L_{\phi}$ to $\bP(W)$ is the union of 
the isomorphisms $\phi'$ and $\phi''$, say, $\phi=\phi'\cup\phi''$. Let
$\psi_1=2s$ and $\psi_2=2\nu s+t$. Then $\phi'$ is induced by the homomorphism 
$(\phi')^*\in\Hom(W,U)$ such that $(\phi')^*(w_j)=\psi_j$. 
On the other hand let $U''_{\phi}=\bC\lambda+\bC t$, 
$\psi''_1=2t$ and $\psi''_2=\lambda+2\nu t$ where we note 
$\psi''_j$ is the linear part of $\phi_j$ in $t$ with $s=1$. 
Then $C''_{\phi}=\bP(U''_{\phi})$ and 
$\phi''$ is induced by the homomorphism
$(\phi')^*\in\Hom(W,U''_{\phi})$ such that $(\phi')^*(w_j)=\psi''_j$.
Furthermore the pull back by $\phi^*$ of the normal sequence for $L$ 
\begin{equation*}
0\to \phi^*N_{L/M}\to \phi^*N_{L/\bP} \overset{\phi^*D_L}\longrightarrow \phi^*O_L(k) \to 0, 
\end{equation*}
yields exact sequences with natural vertical homomorphisms:
\begin{equation*}
\CD
0 @>>> \phi^*N_{L/M} @>>>(\phi')^*N_{L/M}\oplus (\phi'')^*N_{L/M}@>>> \bC @>>> 0\\
@VV{}V @VV{}V @VVV @VV{}V@VVV \\
0 @>>> \phi^*N_{L/\bP} @>>>(\phi')^*N_{L/\bP}\oplus (\phi'')^*N_{L/\bP}
@>>> V^{\vee}/W^{\vee} @>>> 0\\
@VV{}V @VV{}V @VVV @VV{}V@VVV \\
0 @>>> \phi^*O_L(k) @>>>O_{C'_{\phi}}(k)\oplus O_{C''_{\phi}}(k)@>>> \bC @>>> 0.
\endCD
\end{equation*}

This yields the following long exact sequences:
\begin{equation*}
\CD
0 @>>> H^0((\phi')^*N_{L/M}) @>>>U\otimes V^{\vee}/W^{\vee}@>H^0((\phi')^*D_L)>> S^kU  \\
@>>> H^1((\phi')^*N_{L/M})@>>> 0\\
0 @>>> H^0((\phi'')^*N_{L/M}) @>>>U''_{\phi}\otimes V^{\vee}/W^{\vee}@>H^0((\phi'')^*D_L)>> S^kU''_{\phi}\\ 
@>>>H^1((\phi'')^*N_{L/M})@>>> 0\\
\endCD
\end{equation*}
whence $H^1((\phi')^*N_{L/M})=H^1((\phi'')^*N_{L/M})=0$, and both $H^0((\phi')^*N_{L/M})$ 
and $H^0((\phi'')^*N_{L/M})$ are one-dimensional. Let $U'$ be the subspace of $U$ 
consisting of elements vanishing at $C'_{\phi}\cap C''_{\phi}$, namely 
the subspace spanned by $t$. 
Then the restriction of $H^0((\phi')^*D_L)$ to $U'\otimes V^{\vee}/W^{\vee}$
equals $t\cdot H^0((\phi')^*D^{-}_L)$. Hence 
\begin{align*}
\coker H^0(\phi^*D_L)&\simeq t\cdot S^7U/t\cdot\image H^0((\phi')^*D^{-}_L)
\oplus\coker H^0((\phi'')^*D_L)\\
&\simeq S^7U/t\cdot\image H^0((\phi')^*D^{-}_L)
\simeq\coker H^0((\phi')^*D^{-}_L).
\end{align*}

One could understand the above isomorphism as 
\begin{equation*}\coker H^0(\phi^*D_L)
=\coker (\phi)^*H^0(D^{-}_L)\otimes(S^2U/\phi^*W).
\end{equation*}

Thus $H^0(\phi^*N_{L/M})$ is one-dimensional, 
while $H^1(\phi^*N_{L/M})$ is 3-dimensional.
\par

This is immediately generalized into the following
\begin{lemma}\label{lemma:ker and coker}For any $\phi^*\in X^0$ we have
\begin{align*}
\Ker H^0(\phi^*D_L)&=\phi^*\Ker H^0(D_L),\\
\coker H^0(\phi^*D_L)
&=(\phi^*\coker H^0(D^{-}_L))\otimes (S^2U/\phi^*W).
\end{align*}
\end{lemma}
\begin{lemma}We define a line bundle \/ $\bL_0$ (resp. $\bL_1$) on $Y\ (\simeq \bP(S^2W))$ 
by the assignment:
\begin{align*}
X^0\ni\phi^*\mapsto \phi^*\Ker H^0(D_L)\ (resp. \phi^*\coker H^0(D^{-}_L)).
\end{align*}
Then $\bL_k\simeq O_{\bP(S^2W)}$. \end{lemma}
\begin{pf}
We know that $\phi^*\Ker H^0(D_L)$ is generated by $\phi_2e_3^{\vee}-\phi_1e_4^{\vee}$. 
By the $\SL(2)$-variable change of $s$ and $t$, $\phi_j$ 
is transformed into a new quadratic polynomial, which is however the same as the first $\phi_j$.
This shows the generator is unchanged, whence $\bL_0\simeq O_{\bP(S^2W)}$. 
The proof for $\bL_1$ is the same.
\qed
\end{pf}

\begin{lemma}\label{lemma:O(-1/2)}
We define a coherent sheaf \/ $\bL$ on the stack $Y\ (\simeq \bP(S^2W))$ 
(See Remark below)
by the assignment:
\begin{align*}
X^0\ni\phi^*\mapsto S^2U/\phi^*W.
\end{align*}
Then $\bL^2\simeq O_{\bP(S^2W)}(-1)$.
\end{lemma}
\begin{pf}
The GIT-quotient $Y^0$ is covered with the images of $X'_j$:
\begin{align*}
X'_1&=\{(\phi_1,\phi_2)\in X^0;
\phi_1=2st,\ \phi_2=\lambda s^2+2\nu st+t^2,\ 
\lambda,\nu\in\bC\},\\
X'_2&=\{(\phi_1,\phi_2)\in X^0;
\phi_1= p s^2+2q st+t^2,\ \phi_2=2 st,\ p, q\in\bC\}.
\end{align*} 

It is clear that the natural image of $X'_j$ in $Y$ is $Y_j$. 
The map $\phi$ given by $\phi^*=(\phi_1,\phi_2)\in Y_1$ has natural 
$\bZ_{2}$ involution generated by,
\begin{equation*}
r:\;(\sqrt{\lambda}s+t,\sqrt{\lambda}s-t)\rightarrow 
(\sqrt{\lambda}s+t,-(\sqrt{\lambda}s-t)).
\label{ref}
\end{equation*} 
Since 
\begin{eqnarray*}
 2st&=&\frac{1}{2\sqrt{\lambda}}((\sqrt{\lambda}s+t)^2-
(\sqrt{\lambda}s-t)^2),\\ 
\lambda s^2+2\nu st+t^2&=&\frac{\nu}{2\sqrt{\lambda}}((\sqrt{\lambda}s+t)^2-
(\sqrt{\lambda}s-t)^2)+\frac{1}{2}((\sqrt{\lambda}s+t)^2+
(\sqrt{\lambda}s-t)^2),
\end{eqnarray*}
it is 
clear that,
\begin{equation*}
r^{*}(\phi_{1})=\phi_{1},\;\;r^{*}(\phi_{2})=\phi_{2},\;\;
r^{*}(\lambda s^2-t^2)=-(\lambda s^2-t^2).
\end{equation*}
Therefore, we can decompose $S^{2}U$ into $\langle\lambda s^2-t^2\rangle_{\bC}\oplus
\langle\phi_{1},\phi_{2}\rangle_{\bC}$ with respect to eigenvalue of $r^{*}$ 
and take $
\lambda s^2-t^2$ as canonical generator of $S^{2}U/\phi^{*}W$. 
 Similarly 
$S^2U/\phi^*W$ is generated by $p s^2-t^2$ on $Y_2$. 
The problem is therefore to write 
$\lambda s^2-t^2$ as an $\Gamma(O_{Y_1\cap Y_2})$-multiple of $pu^2-v^2$ 
when we write
$\phi_2=2uv$ by a variable change in $\GL(2)$. 
The following variable change $(s,t)\mapsto (u,v)$ is in $\GL(2)$:
\begin{equation*}
s=\frac{\sqrt{2\alpha}}{(\beta-\alpha)^2}(2u-\frac{(\beta-\alpha)^2}{2\alpha}
v),\;\; 
t=\frac{\sqrt{2\alpha}}{(\beta-\alpha)^2}(2\beta u-\frac{(\beta-\alpha)^2}{2}
v),
\end{equation*}
where $\alpha,\;\beta$ are roots of the equation $\lambda s^2+2\nu st+t^2=0$.
Under this coordinate change, $\phi_{1}$ and $\phi_{2}$ is rewritten as 
follows:
\begin{equation*}
\phi_{1}=\frac{\lambda}{(\nu^{2}-\lambda)^{2}}u^2+2\frac{\nu}
{\nu^{2}-\lambda}uv+v^2=:
pu^{2}+2quv+v^2, \; \; \phi_2=2uv.
\end{equation*}
Then we have 
\begin{align*}
pu^{2}-v^{2}=-\frac{2}{\beta-\alpha}(\lambda s^{2}-t^{2})
=-\frac{1}{\sqrt{\nu^{2}-\lambda}}(\lambda s^{2}-t^{2})
=-\sqrt{\frac{D_1}{D_2}}(\lambda s^{2}-t^{2}).
\end{align*}

Similarly by computing the effect on $S^2U/\phi^*W$ by the variable change 
from $X'_1$ into $X'_0$, we see that $\bL^2$ is isomorphic to 
$O_{\bP(S^2W)}(-1)$. This completes the proof. 
\qed
\end{pf}

\begin{rem}
We remark that the space $X$ must be regared as a $\bQ$-stack $Y^{stack}$ as follows:
First we define $\phi_0=-\phi_1-\phi_2$. For each atlas $X'_\alpha$ we define 
an atlas $Y^{stack}_\alpha$ $(\alpha=0,1,2)$ by
\begin{eqnarray*}
Y^{stack}_0&=&\{(\phi_0,\phi_1,\phi_2,\pm\psi_0)\in X^{0}\times S^2U;
\phi_0=2st,\ \phi_1=a s^2+2b st+t^2,\\
&&\psi_0=a s^2-t^2\ a, b\in\bC\},\\
Y^{stack}_1&=&\{(\phi_0,\phi_1,\phi_2,\pm\psi_1)\in X^{0}\times S^2U;
\phi_1=2st,\ \phi_2=\lambda s^2+2\nu st+t^2,\\
&&\psi_1=\lambda s^2-t^2\ 
\lambda,\nu\in\bC\},\\
Y^{stack}_2&=&\{(\phi_0,\phi_1,\phi_2,\pm\psi_2))\in X^{0}\times S^2U;
\phi_1= p s^2+2q st+t^2,\ \phi_2=2 st,\\
&&\psi_2=ps^2-t^2, p, q\in\bC\}.
\end{eqnarray*} 

Since 
$\bL^2\simeq O_{\bP(S^{2}W)}(-1)$
we have $c_1(\bL)=-\frac{1}{2}c_1(O_{\bP(S^{2}W)}(1))$
in the Chow ring $A(Y^{stack})_{\bQ}=A(X)_{\bQ}=A(\bP(S^{2}W))_{\bQ}$.
\end{rem}

\section{Proof of the main theorem}
\begin{thm}
\begin{eqnarray*}
&&\pi_{*}(c_{top}(H^{1}))=\frac{1}{8}\biggl[\frac{c(S^{k-1}Q)}
{1-\frac{1}{2}c_{1}(Q)}
\biggr]_{k-N},
\end{eqnarray*}
where $\pi$ is the natural projection from $\bar{M}_{0,0}(L,2)$ 
to $G$ and $[*]_{k-N}$ is the operation of picking up 
the degree $2(k-N)$ part of Chern classes.  
\end{thm}
\begin{pf} From now on we denote the coherent sheaf $\bL$ 
in Lemma~\ref{lemma:O(-1/2)} by $O_{\bP}(-\frac{1}{2})$.
In view of the results from the previous section, what remains is to evaluate 
the top chern class of $\bigl(S^{k-1}Q/((V^{\vee}\otimes O_{G})/
Q^{\vee})\bigr)\otimes O_{\bP}(-\frac{1}{2})$ on $\bP(S^{2}Q)$.
Since double cover maps parametrized by $\bP(S^{2}Q)$ 
have natural 
$\bZ_{2}$ involution $r$ given in the previous section, 
we have to multiply the result of integration 
on $\bP(S^{2}Q)$ by the factor $\frac{1}{2}$ \cite{thber}, \cite{fp}.
With this set-up, let $\pi^{\prime} :\bP(S^{2}Q)\rightarrow G$ be the natural 
projection.   
Then what we have to compute is $\pi_{*}(c_{top}(H^{1}))=
\frac{1}{2}\pi^{\prime}_{*}(c_{top}(H^{1}))=
\frac{1}{2}\pi^{\prime}_{*}(c_{top}(\bigl(S^{k-1}Q/
((V^{\vee}\otimes O_{G})/
Q^{\vee})\bigr)\otimes O_{\bP}(-\frac{1}{2})))$.  Let $z$ be 
$c_{1}(O_{\bP}(1))$. Then we obtain,
\begin{eqnarray*}
&&\frac{1}{2}\pi^{\prime}_{*}(c_{top}(\bigl(S^{k-1}Q/((V^{\vee}\otimes O_{G})/
Q^{\vee})\bigr)\otimes O_{\bP}(-\frac{1}{2})))\nonumber\\
&&=\frac{1}{2}\sum_{j=0}^{k-N+2}c_{k-N+2-j}(S^{k-1}Q\oplus Q^{\vee})\cdot
\pi_{*}(z^{j})\cdot(-\frac{1}{2})^{j}\nonumber\\
&&=\frac{1}{8}\sum_{j=0}^{k-N}c_{k-N-j}(S^{k-1}Q\oplus Q^{\vee})\cdot
s_{j}(S^{2}Q)\cdot(-\frac{1}{2})^{j}\nonumber\\
&&=\frac{1}{8}\biggl[\frac{c(S^{k-1}Q)\cdot c(Q^{\vee})}
{1-\frac{1}{2}c_{1}(S^{2}Q)+\frac{1}{4}c_{2}(S^{2}Q)-
\frac{1}{8}c_{3}(S^{2}Q)}\biggr]_{k-N},
\end{eqnarray*}
where $s_{j}(S^{2}Q)$ is the $j$-th 
Segre class of $S^{2}Q$. But if we decompose 
$c(Q)$ into $(1+\alpha)(1+\beta)$, we can easily see,
\begin{eqnarray*}
\frac{c(Q^{\vee})}
{1-\frac{1}{2}c_{1}(S^{2}Q)+\frac{1}{4}c_{2}(S^{2}Q)-
\frac{1}{8}c_{3}(S^{2}Q)}&=&\frac{(1-\alpha)(1-\beta)}
{(1-\alpha)(1-\frac{1}{2}(\alpha+\beta))(1-\beta)}\\
&=&\frac{1}{1-\frac{1}{2}c_{1}(Q)}.
\end{eqnarray*}
\hspace{30em}\qed 
\end{pf}

Finally, by combining the above theorem with the divisor axiom of 
Gromov-Witten invariants, we can prove the decomposition formula of 
degree $2$ rational Gromov-Witten invariants of $M_{N}^{k}$ found 
from numerical experiments.
 \begin{cor}
\begin{eqnarray*}
&& \langle{\cal O}_{e^{a}}{\cal O}_{e^{b}}{\cal O}_{e^{c}}\rangle_{0,2}=
\langle{\cal O}_{e^{a}}{\cal O}_{e^{b}}{\cal O}_{e^{c}}\rangle_{0,2
\rightarrow 2}
+8\langle\pi_{*}(c_{top}(H^{1}))
{\cal O}_{e^{a}}{\cal O}_{e^{b}}{\cal O}_{e^c}\rangle_{0,1},
\end{eqnarray*}
where $\langle{\cal O}_{e^{a}}{\cal O}_{e^{b}}{\cal O}_{e^{c}}\rangle_{0,
2\rightarrow 2}$
is the number of conics that intersect cycles Poincar\'e
dual to $e^{a}$, $e^{b}$ 
and $e^{c}$. We also denote by $\langle\pi_{*}(c_{top}(H^{1}))
{\cal O}_{e^{a}}{\cal O}_{e^{b}}{\cal O}_{e^c}\rangle_{0,1}$ 
the integral: 
$$ \int_{G(2,V)} c_{top}(S^{k}Q)\wedge\pi_{*}(c_{top}(H^{1}))
\wedge\sigma_{a-1}\wedge\sigma_{b-1}\wedge\sigma_{c-1}.$$  
\end{cor}
\section{generalization to twisted cubics}

In this section, we present a decomposition formula of degree $3$ 
rational Gromov-Witten invariants found from numerical experiments 
using the results of \cite{ES}. 
\begin{conjecture}
If $k-N=1$, we have the following equality:
\begin{eqnarray*}
&&\pi_{*}(c_{top}(H^{1}))=\\
&&\frac{1}{27}\biggl
(\bigl(\frac{1}{24}(27k^{2}-55k+26)k(k-1)+\frac{2}{9}\bigr)c_{1}(Q)^{2}+
\bigl(\frac{7}{6}(k+1)k(k-1)+\frac{1}{9}\bigr)c_{2}(Q)\biggr).
\end{eqnarray*}
where $\pi: \overline{M}_{0,0}(L,3)\rightarrow 
\overline{M}_{0,0}(M_{N}^{k},1)$ is the natural projection.
\end{conjecture}
In the $k-N>1$ case, we 
have not found the explicit formula, because in the $d=3$ case, 
we have another contribution from multiple cover maps of type 
$(2+1)\rightarrow (1+1)$. Here multiple cover map of type 
$(2+1)\rightarrow (1+1)$ is the map from nodal curve $\bP^{1}\vee\bP^{1}$ 
to nodal conic $L_{1}\vee L_{2} \subset M_{N}^{k}$, that maps the first 
(resp. the second) $\bP^{1}$ to $L_{1}$ (resp. $L_{2}$) by two to one (resp. 
one to one). 
In the $k-N=1$ case, we have also determined the contributions from multiple 
cover maps of $(2+1) \to (1+1)$ to nodal conics. 
\begin{cor}
If $k-N=1$, 
$\langle{\cal O}_{e^{a}}{\cal O}_{e^{b}}{\cal O}_{e^{c}}\rangle_{0,3}$
is decomposed into the following contributions: 
\begin{eqnarray*}
&& \langle{\cal O}_{e^{a}}{\cal O}_{e^{b}}{\cal O}_{e^{c}}\rangle_{0,3}=
\langle{\cal O}_{e^{a}}{\cal O}_{e^{b}}{\cal O}_{e^{c}}\rangle_{0,3
\rightarrow3}
+\frac{1}{k}\bigl(
\frac{9}{4}\langle{\cal O}_{e^{a}}{\cal O}_{e^{b}}{\cal O}_{e^{c}}
{\cal O}_{e^{3}}\rangle_{0,1}\langle{\cal O}_{e^{N-5}}\rangle_{0,1}\\
&&+\frac{3}{2}\langle{\cal O}_{e^{a}}{\cal O}_{e^{b}}
{\cal O}_{e^{c+2}}\rangle_{0,1}
\langle{\cal O}_{e^{N-c-4}}{\cal O}_{e^{c}}\rangle_{0,1}+
\frac{3}{2}\langle{\cal O}_{e^{b}}{\cal O}_{e^{c}}
{\cal O}_{e^{a+2}}\rangle_{0,1}
\langle{\cal O}_{e^{N-a-4}}{\cal O}_{e^{a}}\rangle_{0,1}\\
&&+\frac{3}{2}\langle{\cal O}_{e^{c}}{\cal O}_{e^{a}}
{\cal O}_{e^{b+2}}\rangle_{0,1}
\langle{\cal O}_{e^{N-b-4}}{\cal O}_{e^{b}}\rangle_{0,1}\bigr)
\\
&&+27\langle\pi_{*}(c_{top}(H^{1}))
{\cal O}_{e^{a}}{\cal O}_{e^{b}}{\cal O}_{e^c}\rangle_{0,1},
\end{eqnarray*}
where $\langle{\cal O}_{e^{a}}{\cal O}_{e^{b}}{\cal O}_{e^{c}}\rangle_{0,
3\rightarrow 3}$
is the number of twisted cubics that intersect cycles Poincar\'e
dual to $e^{a}$, $e^{b}$ and $e^{c}$. 
\end{cor}
\begin{pf}
In the $k-N=1$ case, dimension of moduli space of multiple 
cover maps of $(2+1) \to (1+1)$ to nodal conics is given by 
$N-6+N-6-(N-4)+2=N-6$, hence the rank of $H^{1}$ is given by $N-6-(N-5-3)=2$. 
On the other hand, dimension of moduli space of $d=2$ multiple cover maps
of $\bP^{1}\to\bP^{1}$ is $2$, the degree of the form of 
$\tilde{\pi}_{*}(c_{top}(H^{1}))$ 
equals to $2-2=0$, where $\tilde{\pi}$ is the projection map 
that projects out the fiber locally isomorphic to the moduli space of 
$d=2$ multiple cover maps. This situation is exactly the same as the 
Calabi-Yau case. Therefore, we can use the well-known result by Aspinwall 
and Morrison, that says for $n$-point rational Gromov-Witten invariants 
for Calabi-Yau manifold, $\tilde{\pi}_{*}(c_{top}(H^{1}))$ 
for degree $d$ multiple 
cover map is given by,
$$
\tilde{\pi}_{*}(c_{top}(H^{1}))=\frac{1}{d^{3-n}}.
$$
With this formula, we add up all the combinatorial possibility of 
insertion of external operator ${\cal O}_{e^{a}}$, ${\cal O}_{e^{b}}$ and
${\cal O}_{e^{c}}$,
\begin{eqnarray*}
&&\frac{1}{k}\bigl(\langle\tilde{\pi}_{*}(c_{top}(H^{1}))
{\cal O}_{e^{a}}{\cal O}_{e^{b}}{\cal O}_{e^{c}}
{\cal O}_{e^{3}}\rangle_{0,1}\langle{\cal O}_{e^{N-5}}\rangle_{0,1}\\
&&+\langle\tilde{\pi}_{*}(c_{top}(H^{1})){\cal O}_{e^{a}}{\cal O}_{e^{b}}
{\cal O}_{e^{c+2}}\rangle_{0,1}
\langle{\cal O}_{e^{N-c-4}}{\cal O}_{e^{c}}\rangle_{0,1}\nonumber\\
&&+\langle\tilde{\pi}_{*}(c_{top}(H^{1})){\cal O}_{e^{b}}{\cal O}_{e^{c}}
{\cal O}_{e^{a+2}}\rangle_{0,1}
\langle{\cal O}_{e^{N-a-4}}{\cal O}_{e^{a}}\rangle_{0,1}\\
&&+\langle\tilde{\pi}_{*}(c_{top}(H^{1})){\cal O}_{e^{c}}{\cal O}_{e^{a}}
{\cal O}_{e^{b+2}}\rangle_{0,1}
\langle{\cal O}_{e^{N-b-4}}{\cal O}_{e^{b}}\rangle_{0,1}\\
&&+\langle{\cal O}_{e^{a}}{\cal O}_{e^{b}}
{\cal O}_{e^{c+2}}\rangle_{0,1}
\langle{\cal O}_{e^{N-c-4}}\tilde{\pi}_{*}(c_{top}(H^{1})){\cal O}_{e^{c}}
\rangle_{0,1}\nonumber\\
&&+\langle{\cal O}_{e^{b}}{\cal O}_{e^{c}}
{\cal O}_{e^{a+2}}\rangle_{0,1}
\langle{\cal O}_{e^{N-a-4}}\tilde{\pi}_{*}(c_{top}(H^{1})){\cal O}_{e^{a}}
\rangle_{0,1}\\
&&+\langle{\cal O}_{e^{c}}{\cal O}_{e^{a}}
{\cal O}_{e^{b+2}}\rangle_{0,1}
\langle{\cal O}_{e^{N-b-4}}\tilde{\pi}_{*}(c_{top}(H^{1})){\cal O}_{e^{b}}
\rangle_{0,1}\\
&&+\langle{\cal O}_{e^{a}}{\cal O}_{e^{b}}{\cal O}_{e^{c}}
{\cal O}_{e^{3}}\rangle_{0,1}\langle
{\cal O}_{e^{N-5}}\tilde{\pi}_{*}(c_{top}(H^{1}))\rangle_{0,1}\bigr)\\
&&=\frac{1}{k}\bigl(2\langle
{\cal O}_{e^{a}}{\cal O}_{e^{b}}{\cal O}_{e^{c}}
{\cal O}_{e^{3}}\rangle_{0,1}\langle{\cal O}_{e^{N-5}}\rangle_{0,1}\\
&&
+\langle{\cal O}_{e^{a}}{\cal O}_{e^{b}}
{\cal O}_{e^{c+2}}\rangle_{0,1}
\langle{\cal O}_{e^{N-c-4}}{\cal O}_{e^{c}}\rangle_{0,1}
+\langle{\cal O}_{e^{b}}{\cal O}_{e^{c}}
{\cal O}_{e^{a+2}}\rangle_{0,1}
\langle{\cal O}_{e^{N-a-4}}{\cal O}_{e^{a}}\rangle_{0,1}\\
&&+\langle{\cal O}_{e^{c}}{\cal O}_{e^{a}}
{\cal O}_{e^{b+2}}\rangle_{0,1}
\langle{\cal O}_{e^{N-b-4}}{\cal O}_{e^{b}}\rangle_{0,1}\\
&&+\frac{1}{2}\langle
{\cal O}_{e^{a}}{\cal O}_{e^{b}}{\cal O}_{e^{c+2}}
\rangle_{0,1}\langle
{\cal O}_{e^{N-c-4}}{\cal O}_{e^{c}}\rangle_{0,1}
+\frac{1}{2}\langle{\cal O}_{e^{b}}{\cal O}_{e^{c}}
{\cal O}_{e^{a+2}}\rangle_{0,1}
\langle{\cal O}_{e^{N-a-4}}\pi_{*}{\cal O}_{e^{a}}\rangle_{0,1}\\
&&+\frac{1}{2}\langle{\cal O}_{e^{c}}{\cal O}_{e^{a}}
{\cal O}_{e^{b+2}}\rangle_{0,1}
\langle{\cal O}_{e^{N-b-4}}{\cal O}_{e^{b}}\rangle_{0,1}\\
&&+\frac{1}{4}\langle{\cal O}_{e^{a}}{\cal O}_{e^{b}}{\cal O}_{e^{c}}
{\cal O}_{e^{3}}\rangle_{0,1}\langle
{\cal O}_{e^{N-5}}\rangle_{0,1}\bigr).\\
\end{eqnarray*}
The last expression is nothing but the formula we want.
\qed
\end{pf}
\\
\\
{Masao Jinzenji${}^{\dagger}$, Iku Nakamura${}^{\star}$, 
Yasuki Suzuki}\\
{\it Division of Mathematics, Graduate School of Science,
 Hokkaido University \\
\it  Kita-ku, Sapporo, 060-0810, Japan\\
{\it e-mail address: $\dagger$ jin@math.sci.hokudai.ac.jp, $\star$
nakamura@math.sci.hokudai.ac.jp}
\\
\\

\end{document}